\documentclass[12pt]{amsart}
\usepackage{amssymb}
\usepackage{epsf}
\usepackage{colordvi}
\usepackage{multirow}
\usepackage{graphicx}
\usepackage{amssymb}

\newcommand{\bC}{\mathbb{C}}
\renewcommand{\P}{\mathbb{P}}
\newcommand{\QED}{\hfill\raisebox{-5pt}{\includegraphics[height=15.5pt]{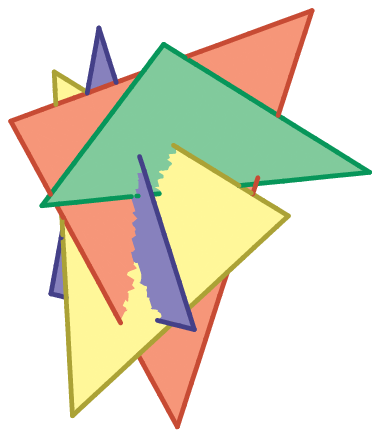}}\medskip}

\newcounter{FNC}[page]
\def\newfootnote#1{{\addtocounter{FNC}{2}$^\fnsymbol{FNC}$%
     \let\thefootnote\relax\footnotetext{$^\fnsymbol{FNC}$#1}}}

\newtheorem{theorem}{Theorem}
\newtheorem{lemma}[theorem]{Lemma}

\newtheorem{rem}[theorem]{Remark}
\newenvironment{remark}{\begin{rem}\rm}{\end{rem}}

\copyrightinfo{}{}

\title[Line tangents to four triangles in three-dimensional space]{Line tangents
  to four triangles\\ in three-dimensional space} 

\author[Br{\"o}nnimann]{Herv{\'e}  Br{\"o}nnimann}
\address{CIS Dept\\
         Polytechnic University\\
         Six Metrotech\\
         Brooklyn  NY 11201\\
         USA}
\email{{\tt hbr@poly.edu}}

\author[Devillers]{Olivier Devillers}
\address{INRIA\\
         BP 93\\
         06902 Sophia-Antipolis\\
         France}
\email{{\tt olivier.devillers@sophia.inria.fr}}

\author[Lazard]{Sylvain Lazard}
\address{INRIA - LORIA  \\
         54506  Vand{\oe}uvre-l{\`e}s-Nancy Cedex\\
         France}
\email{{\tt sylvain.lazard@loria.fr}}

\author[Sottile]{Frank Sottile}
\address{Department of Mathematics\\
         Texas A\&M University\\
         College Station, TX 77843\\
         USA}
\email{sottile@math.tamu.edu}

\thanks{Research of Br{\"o}nnimann supported in part by NSF CAREER Grant CCR-0133599}
\thanks{Research of Sottile is supported in part by NSF CAREER grant DMS-0134860}
\thanks{Research initiated at the Second McGill-INRIA Workshop on Computational
  Geometry in Computer Graphics, February 7--14 2003 co-organized by
  H. Everett, S. Lazard, and S. Whitesides, and held at the Bellairs
  Research Institute of McGill University in Holetown, St. James,
  Barbados, West Indies. A preliminary version appeared in Proc. Canad. Conf. Comput. Geom., Montreal, August 2004.}

\begin{document}


\begin{abstract}
We investigate the lines tangent to four triangles in $\mathbb{R}^3$.
By a construction, there can be as many as 62 tangents.
We show that there are at most 162 connected components of tangents, and at most
156 if the triangles are disjoint.
In addition, if the triangles are in (algebraic) general position, then the number of
tangents is finite and it is always even.
\end{abstract}
\maketitle

\section*{Introduction}

Motivated by visibility problems, we investigate lines tangent to four
triangles in $\mathbb{R}^3$.  In computer graphics and robotics, scenes
are often represented as unions of not necessarily disjoint polygonal or
polyhedral objects.  The objects that can be seen in a particular
direction from a moving viewpoint may change when the line of sight
becomes tangent to one or more objects in the scene. Since this line of
sight is tangent to a subset of the edges of the polygons and polyhedra
representing the scene, we are also led to questions about lines tangent
to segments and to polygons.  Four polygons will typically have finitely
many common tangents, while 5 or more will have none and 3 or fewer will
have either none or infinitely many.

This paper is the third in a series of papers by the authors and their
collaborators investigating such questions.  The paper~\cite{socg04}
investigated the lines of sight tangent to four convex polyhedra in a
scene of $k$ convex but not necessarily disjoint polyhedral objects, and
proved that there could be up to but no more than $\Theta(n^2k^2)$
connected components of such lines. The same bound for the
considerably easier case of disjoint convex polyhedra in algebraic
general position was proved earlier~\cite{Efrat, cccg02}.  We would
like, however, 
to investigate how high the constants hidden in the $O()$ notation are.
The paper~\cite{BELSW} offers a detailed study of transversals to $n$
line segments in $\mathbb{R}^3$ and proved that although there are at
most 2 such transversals for four segments in (algebraic) general
position, there are at most $n$ such connected components of
transversals in any case. In this paper, we consider the 
case of four triangles in
$\mathbb{R}^3$, and establish lower and upper bounds on the number of
tangent lines.

A \Blue{\emph{triangle}} in $\mathbb{R}^3$ is the convex hull of three
distinct (and non-collinear) points in $\mathbb{R}^3$.  A line is
\Blue{\emph{tangent}} to a triangle if it meets an edge of the triangle.
Note that a line tangent to each of four triangles forming a scene
corresponds to an unoccluded line of sight in that scene. If there are
$k>4$ triangles, then the bound $\Theta(k^4)$ of~\cite{socg04} stands
(as the total number of edges is $n=3k$ and one of the lower bound
examples is made of triangles). We thus investigate the case of four
triangles. Let $n(t_1,t_2,t_3,t_4)$ be number of lines tangent to  four
triangles $t_1$, $t_2$, $t_3$, and $t_4$ in $\mathbb{R}^3$.  This number
may be infinite if the lines supporting the edges of the different
triangles are not in general position.

Our first step is to consider the algebraic relaxation of 
this geometric problem in which we replace each edge of a triangle by the line
in $\bC\P^3$ supporting it, and then ask for the set of lines in  $\bC\P^3$
which meet one supporting line from each triangle.
Since there are $3^4=81$ such quadruples of supporting lines, this 
is the disjunction of 81 instances of the classical problem of transversals
to four given lines in $\bC\P^3$.
As there are two such transversals to four given lines in general position,
\begin{figure}[htb]
 \begin{center}
   \setlength{\unitlength}{1pt}
   \begin{picture}(270,160)(-5,-5)
    \put(  5,  5){\includegraphics[height=2in]{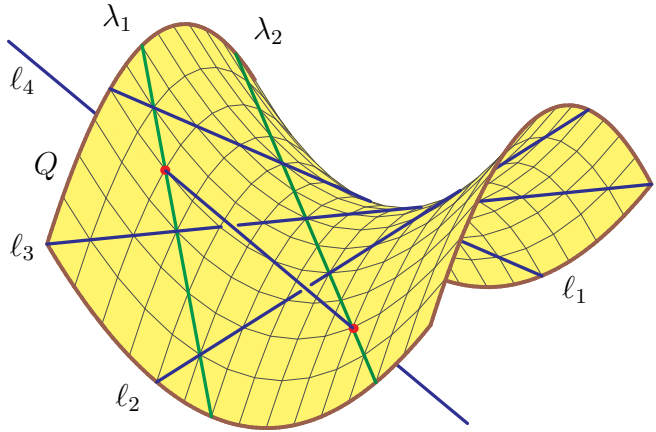}} 
    \put( 35,  8){$\ell_2$}
    \put( -5, 65){$\ell_3$}
    \put(204, 49){$\ell_1$}
    \put( -5,128){$\ell_4$}
    \put( 30,152){$\lambda_1$}
    \put( 87,147){$\lambda_2$}
    \put(  4,96){$Q$}
   \end{picture}
 \caption{The lines $\ell_1$, $\ell_2$ and $\ell_3$ span a hyperbolic paraboloid $Q$
    which meets line $\ell_4$ in two points.  The two lines $\lambda_1$ and $\lambda_2$
    are the transversals to the four lines $\ell_1$ $\ell_2$, $\ell_3$, and~$\ell_4$.}
 \label{fig:one}
 \end{center}
\end{figure}
we expect that this algebraic relaxation has 162 solutions.
We say that four triangles $t_1,t_2,t_3,t_4$ are in (algebraic) 
\Blue{\emph{general position}}
if each of the 81 quadruples of supporting lines have two transversals in
$\bC\P^3$ and all 162
transversals are distinct.  
Let $\mathcal{T}$ be the configuration space of all quadruples of triangles in
$\mathbb{R}^3$ and $T\subset\mathcal{T}$ consist of those quadruples which are in
general position.
Thus if $(t_1,t_2,t_3,t_4)\in T$,  the number $n(t_1,t_2,t_3,t_4)$ is finite and
is at most $162$.

Our first result is a congruence.

\begin{theorem}\label{T:even}
 If $(t_1,t_2,t_3,t_4)\in T$, then $n(t_1,t_2,t_3,t_4)$ is even.
\end{theorem}

Our primary interest is the number
 \[ 
    N\ :=\  \max\{ n(t_1,t_2,t_3,t_4)\mid (t_1,t_2,t_3,t_4)\in T\}\,.
 \]
Our results about this number $N$ are two-fold. First, we show that $N\geqslant 62$.

\begin{theorem}\label{T:lower}
  There are four disjoint triangles in $T$ with $62$ common tangent lines.  
\end{theorem}

The idea is to perturb a configuration of four lines in $\mathbb{R}^3$ with two real
transversals, such as in Figure~\ref{fig:one}. 
The triangles in our construction are very `thin'---the smallest angle
among them measures about $10^{-11}$ degrees.
We ran a computer search for `fatter' triangles having many common tangents,
checking the number of tangents to 5 million different quadruples of triangles.
Several had as many as 40 common tangents.
This is discussed in Section~\ref{S:search}.

We can improve the upper bound on $N$ when the triangles are
disjoint.

\begin{theorem}\label{T:upper}
  Four triangles in $T$ admit at most $162$ distinct common tangent lines.
  This number is at most $156$ if the triangles are disjoint.
\end{theorem}

When the four triangles are not in general position, the 
number of tangent lines can be infinite. In this case, we
may group these tangents by connected components: two line tangents are in
the same component if one may move continuously between the two lines
while staying tangent to the four triangles.
Each quadruple of edges may induce up to four components of tangent
lines~\cite{BELSW}, giving a trivial upper bound of $324$. 
This may be improved.

\begin{theorem}\label{T:upper-comp}
  Four triangles have at most $162$ connected components of common tangents.
  If the triangles are disjoint, then this number is at most $156$.
\end{theorem}

We believe that these upper bounds are far from optimal.
Theorems~\ref{T:even},~\ref{T:lower},~\ref{T:upper}, and~\ref{T:upper-comp} are
proved in Sections~\ref{sec:even},~\ref{sec:lower},~\ref{sec:upper},
and~\ref{S:upper-comp}, respectively.
Section~\ref{S:search} discusses our search for `fat' triangles with many common  
tangents.


\section{A congruence}\label{sec:even}

We prove Theorem~\ref{T:even} by showing that any two quadruples of triangles
whose supporting lines are in general position are connected by a path such that
common tangents are created and destroyed in pairs along that path. 
Thus the parity of $n(t_1,t_2,t_3,t_4)$ is constant for $(t_1,t_2,t_3,t_4)\in T$.
The theorem follows as there are triangles in $T$ with no common tangents.

We study the complement $\Sigma$ of $T$ in the set $\mathcal{T}$ of quadruples of all
triangles. 
The reason is that the number  $n(t_1,t_2,t_3,t_4)$ of common tangents is constant 
in each connected component of $T$ and so we must pass
through $\Sigma$ to connect quadruples in $T$.
This is called `crossing a wall'.
Since the set of smooth points of $\Sigma$ is open and dense in $\Sigma$, a path may be
found which meets $\Sigma$ only in its smooth points.
Since a smooth point lies on a unique algebraic component, each
crossing involves only one algebraic component of $\Sigma$, and so 
we must also describe the different algebraic components of $\Sigma$.
This amounts to describing what happens near a smooth point of $\Sigma$.
Note that $\Sigma$ is also called the discriminant hypersurface of $\mathcal{T}$.

Recall that a quadruple $(t_1,t_2,t_3,t_4)$ lies in $T$ only if
 \begin{enumerate}
  \item[(A)] There are two lines in $\bC\P^3$ transversal to each quadruple
             $\ell_1,\ell_2,\ell_3,\ell_4$ of lines supporting one edge from each
             triangle, and 
  \item[(B)] the 162 such lines are distinct.
 \end{enumerate}

\begin{lemma}\label{L:discriminant}
Each complex algebraic component of $\Sigma$ contains an open dense
 set on which exactly one of $(a)$ or $(b)(i)$ or $(b)(ii)$ occurs.
\begin{enumerate}
 \item[$(a)$] There is a unique transversal $\lambda$ in $\bC\P^3$ to one
             quadruple of supporting lines.
 \item[$(b)$] One of the lines $\lambda$ meeting one quadruple of supporting lines
             $\ell_1,\ell_2,\ell_3,\ell_4$ meets one other supporting line 
             $\ell'$. 
             There are two ways for this to occur. Either 
   \begin{enumerate}
     \item[$(i)$] $\lambda$ meets a vertex of the triangle having $\ell'$ as a
               supporting line, or
     \item[$(ii)$] $\lambda$ lies in the plane of the triangle having $\ell'$ as a
               supporting line. 
   \end{enumerate}

\end{enumerate}
 In each case, the distinguished line $\lambda$ is real.
\end{lemma}

\noindent{\it Proof.}
 We consider what happens when one of the conditions (A) or (B) fails, but the
 rest of the configuration remains generic.
 For (A), if there is a quadruple $\ell_1,\ell_2,\ell_3,\ell_4$ of supporting lines
 without two common transversals, then either there is only one transversal or there are
 infinitely  many. 
 Since we are considering generic such configurations, we may assume that 
 $\ell_1,\ell_2$, and $\ell_3$ are in general position in that they span a quadric $Q$
 as in Figure~\ref{fig:one}, and ask what happens as $\ell_4$ moves out of
 general position. 
 If $\ell_4$ meets one of $\ell_1,\ell_2$, or $\ell_3$, there still will be two
 lines,  but if $\ell_4$ becomes tangent to $Q$, then there will only be one, as
 the two lines $\lambda_1$ and $\lambda_2$ coalesce.
 Further degeneration is required for there to be infinitely many lines, since
 $\ell_4$ has then to become contained in  $Q$. Thus 
 $(a)$ describes what happens generically when (A) fails for a single quadruple of
 supporting lines.

 For (B), we may assume that each quadruple of supporting lines has two transversals,
 but there are two quadruples with a common transversal.
 The generic way for this to occur is described in $(b)$.
 That is, $\ell_1,\ell_2$, and $\ell_3$ are in general position and a line
 $\lambda$ meeting all three also meets both $\ell_4$ and $\ell'$, and 
 also $\ell_4$ and $\ell'$ are lines supporting edges from the same triangle,
 $t_4$.
 Since $\ell_4$ and $\ell'$ meet in a vertex $v$ of $t_4$ and also span the
 plane $\pi_4$ of $t_4$, either $\lambda$ meets $v$ or else $\lambda$
 lies in $\pi_4$, and these are the two cases $(b)(i)$ and $(b)(ii)$.
 As we consider configurations which are otherwise general, exactly one of these
 two possibilities occurs.

 To see that it is possible for exactly one of $(a)$ or $(b)(i)$ or $(b)(ii)$ to
 occur, begin with a 
 configuration of four triangles in $T$, and allow exactly one supporting
 line of one  triangle to rotate about one vertex, remaining in the plane of the
 triangle. 
 Perturbing the plane of this last triangle if necessary, we see that only
 the configurations described in $(a)$ or $(b)(i)$ or $(b)(ii)$ can occur, each will
 occur finitely many times, and they will occur for distinct angles of rotation.
 This shows that each different possibility
 describes different algebraic components of the discriminant, and that
 each component has an open dense set in which exactly one of these
 possibilities occurs.

 Since the lines and vertices defining the special line $\lambda$
 are all real and $\lambda$ is unique, it will also be real.
\QED

\noindent{\it Proof of Theorem$~\ref{T:even}$.}
 Suppose now that we have two quadruples of triangles in $T$.
 A consequence of Lemma~\ref{L:discriminant} is that there exists a path $\gamma$ in
 $\mathcal{T}$ connecting them such that each time $\gamma$ meets the discriminant
 hypersurface $\Sigma$, exactly one of 
$(a)$ or $(b)(i)$ or $(b)(ii)$ occurs. 
 We need only show that the parity of the number of tangents does not change as we move
 along $\gamma$ and one of $(a)$ or $(b)(i)$ or $(b)(ii)$ occurs.

 If $(a)$ occurs, the number of tangents changes only if the double line $\lambda$
 is tangent to the triangles. 
 Approaching this configuration along the curve $\gamma$, either two real lines or two
 complex lines coalesce into $\lambda$.
 Thus the parity of $n(t_1,t_2,t_3,t_4)$ does not change when crossing $\Sigma$ in a
 component of type $(a)$. 
 
\begin{figure}[htb]
\[
  \setlength{\unitlength}{1.2pt}
   \begin{picture}(120,100)(-5,-7)
    \put(0,0){\includegraphics[height=1.6in]{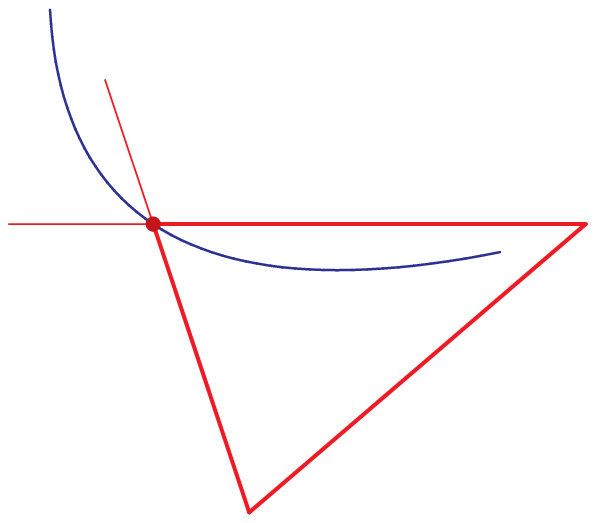}}
    \put(9,88){\vector(4,1){15}}
    \put(-5,75){$C_4$}
    \put(20,47){$v$}
    \put(27,15){$\ell_4$}    \put(65,61){$\ell'_4$}
    \put(75,15){$t_4$}
    \put(47.5,-12){(1)}
   \end{picture}
  \qquad 
  \begin{picture}(120,100)(0,-7)
    \put(0,0){\includegraphics[height=1.6in]{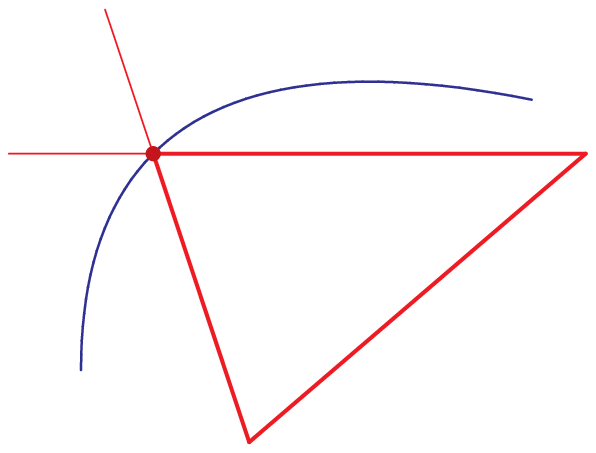}}
    \put(17,23){\vector(4,-1){15}}
    \put(6,40){$C_4$}
    \put(22,67){$v$}
    \put(32,31){$\ell_4$}    \put(72,54){$\ell'_4$}
    \put(84,15){$t_4$}
    \put(47.5,-12){(2)}
  \end{picture}
  \setlength{\unitlength}{1pt}
\]
\caption{Configuration in plane $\pi_4$}
 \label{F:C_4}
\end{figure}

 For $(b)(i)$, we suppose that $\ell'=\ell'_4$ is the line supporting an edge of
 the fourth triangle, $t_4$.
 Let $C_4$ be the conic which is the intersection of the hyperboloid spanned by
 $\ell_1,\ell_2$, and $\ell_3$ with the plane $\pi_4$ spanned by $t_4$.
 Through every point of $C_4$ there is a unique line meeting $\ell_1,\ell_2$, and
 $\ell_3$.
 In particular, the line $\lambda$ corresponds to the vertex $v$ of $t_4$ where
 $\ell_4$ meets $\ell'_4$.
 Figure~\ref{F:C_4} illustrates the two possibilities for the configuration of
 $C_4$ and 
 $t_4$:
 Either (1) $C_4$ meets the interior of $t_4$ or (2) it does not.
 Moving along the curve $\gamma$ perturbs the configuration.
 Topologically, this corresponds to moving $C_4$ off the vertex $v$,
 which is suggested by the arrows in Figure~\ref{F:C_4}.
 In (1), there will be one line near to $\lambda$ meeting $\ell_1,\ell_2$, and $\ell_3$,
 and $t_4$ both before and after the conic  $C_4$ meets the vertex $v$ (but these lines
 will meet different edges of $t_4$).
 In case (2), two lines which meet the supporting lines outside of $t_4$ coalesce into
 $\lambda$, and then become two lines meeting $t_4$.

 For $(b)(ii)$, we also suppose that $\ell'=\ell'_4$ is the line supporting an edge of
 the fourth triangle, $t_4$, and that $\pi_4$ is the plane of the triangle $t_4$.
 Then the three points $\ell_i\cap\pi_4$ for $i=1,2,3$ are collinear and span
 the line $\lambda$ which meets $t_4$.
 We may assume that $\lambda$ does not meet any vertex of the triangle $t_4$.
 But then $\lambda$ meets the interiors of the edges of $t_4$ supported by 
 $\ell_4$ and $\ell'_4$, but not the third edge.
 If we perturb $\ell_1,\ell_2$, and $\ell_3$ to $\ell'_1, \ell'_2$, and
 $\ell'_3$, then there is a line $\mu$ meeting this perturbed triple and
 $\ell_4$ which is close to $\lambda$.
 Similarly, there is a line $\mu'$ meeting this perturbed triple and $\ell'_4$
 which is close to $\lambda$.
 If the points of intersection of $\ell'_1, \ell'_2$, and $\ell'_3$ with $\pi_4$
 are not collinear, then $\mu\neq\mu'$.
 Thus there is no change in the number of common tangents to the four triangles
 when crossing a wall of this type.

 Thus the parity of the number of lines tangent to the four triangles does not change
 when crossing $\Sigma$, which completes the proof of Theorem~\ref{T:even}.
\QED

\section{A construction with 62 tangents}\label{sec:lower}

Consider the four triangles whose vertices are given in Table~\ref{TA:triangles}.\medskip
\begin{table}[htb]

  \begin{center}
   \begin{tabular}{|l||l|}\hline
   \multirow{3}{10pt}{$t_1$}&
   {\scriptsize $(-10.5, 1, -10.5)$}\vspace{-4pt}\\
   &{\scriptsize $(.5628568345479573470378601, 1, .5628568345479573470378601)$}\vspace{-4pt}\\
   &{\scriptsize $(.56285683454726874605620706, .99999999999822994290647247, .56285683454726874605620706)$}\\\hline
   \multirow{3}{10pt}{$t_2$}&
    {\scriptsize $(-10.5, -1, 10.5)$}\vspace{-4pt}\\
   &{\scriptsize $(1.394218989475, -1, -1.394218989475)$}\vspace{-4pt}\\
   &{\scriptsize $(1.3942406911811439954597161, -1.0000237884694881275439271, -1.3942406911811439954597161)$}\\\hline
   \multirow{3}{10pt}{$t_3$}&
    {\scriptsize $(-9.5, -9.5, .25)$}\vspace{-4pt}\\
   &{\scriptsize $(.685825, .685825, .25)$}\vspace{-4pt}\\
   &{\scriptsize $(.69121730616063647303519136, .69121730616063647303519136, .26069756890079842876805653)$}\\\hline
   \multirow{3}{10pt}{$t_4$}&
   {\scriptsize $(9.5, 0, 0)$}\vspace{-4pt}\\
   &{\scriptsize $(-.511, 0, 0)$}\vspace{-4pt}\\
   &{\scriptsize $(-1.0873912730501133759642956, 0, -.51645811088049333541289247)$}\\\hline
  \end{tabular}
 \end{center}

 \caption{Four triangles with 62 common tangents}
 \label{TA:triangles}
\end{table}

\noindent{\bf Theorem~\ref{T:lower}$'$.}
{\it 
 There are exactly $62$ lines tangent to the four triangles of Table~$\ref{TA:triangles}$.
}\medskip

This can be verified by a direct computation.
Software is provided on this paper's web page\newfootnote{{\tt
    http://www.math.tamu.edu/\~{}sottile/stories/4triangles/index.html}}.
More illuminating perhaps is our construction.
The idea is to perturb a configuration of four lines in $\mathbb{R}^3$ with two 
transversals such as in Figure~\ref{fig:one}. 
The resulting triangles of Theorem~\ref{T:lower}$'$ are very thin.
In degrees, their smallest angles are
\[
   t_1\colon 6.482 \times 10^{-12},\ \quad 
   t_2\colon 8.103 \times 10^{-5},\ \quad 
   t_2\colon 4.253 \times 10^{-2},\quad\mbox{and}\quad 
   t_4\colon 2.793\,.
\]

\subsection{The construction}
The lines given parametrically 
\[
   \ell_1\ \colon\ (t,1,t)\,,\quad
   \ell_2\ \colon\ (t,-1,-t)\,,\quad
   \ell_3\ \colon\ (t,t,{\textstyle \frac{1}{4}})\,,\quad\mbox{and}\quad 
   \ell_4\ \colon\ (t,0,0)\,,
\]
have two transversals
\[
  \lambda_1\ \colon\ ({\textstyle \frac{1}{2}}, 2t,t)
   \quad\mbox{and}\quad
  \lambda_2\ \colon\ (-{\textstyle \frac{1}{2}}, 2t,-t)\,.
\]

For each $i=1,2,3,4$, let $Q_i$ be the hyperboloid spanned by the lines
other than $\ell_i$.
For example, $Q_3$ has equation $z=xy$.
The intersection of $Q_i$ with a plane containing $\ell_i$ will be a conic which
meets $\ell_i$ in two points (corresponding to the common transversals
$\lambda_1$ and $\lambda_2$ at $t=\pm\frac{1}{2}$). 
We choose the plane $\pi_i$ so that these two points lie in the same connected
component of the conic.
Here is one possible choice
\[
   \pi_1\ \colon\ x=z\,,\quad
   \pi_2\ \colon\ x=-z\,,\quad
   \pi_3\ \colon\ x=y\,,\quad\mbox{and}\quad 
   \pi_4\ \colon\ y=0\,.
\]

\begin{figure}[htb]
\[
  \begin{picture}(99,103)
   \put(5,10){\includegraphics[width=3cm]{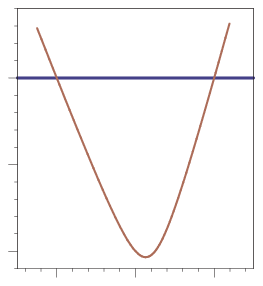}}
   \put(-2,18){$-1$}   \put(-2,45.2){$-\frac{1}{2}$}   \put(6,71){$0$}
   \put(18,1){$-\frac{1}{2}$}   \put(51.75,2){$0$}   \put(75.5,1){$\frac{1}{2}$} 
   \put(49,87){$\pi_1$}
  \end{picture}
     \quad
  \begin{picture}(99,103)
   \put(5,10){\includegraphics[width=3cm]{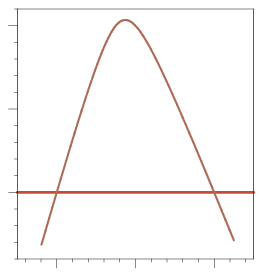}}
   \put(3,34){$0$}   \put(3,61.5){$\frac{1}{2}$}   \put(4,86.5){$1$}
   \put(15.3,1){$-\frac{1}{2}$}   \put(49.75,2){$0$}   \put(75,1){$\frac{1}{2}$}
   \put(77,87){$\pi_2$}
  \end{picture}
     \quad
  \begin{picture}(99,103)
   \put(5,10){\includegraphics[width=3cm]{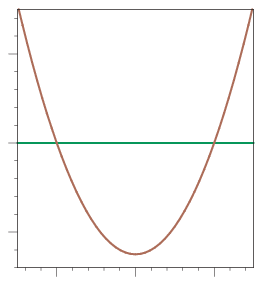}}
   \put(-3,24.5){$-\frac{1}{5}$}   \put(6,51){$0$}  \put(6,79){$\frac{1}{5}$}
   \put(18,1){$-\frac{1}{2}$}   \put(51.75,2){$0$}   \put(75.5,1){$\frac{1}{2}$}
   \put(51,87){$\pi_3$}
  \end{picture}
     \quad
  \begin{picture}(99,103)
   \put(5,10){\includegraphics[width=3cm]{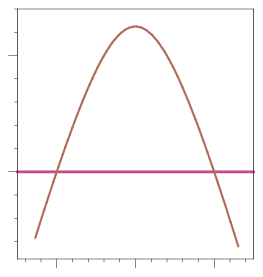}}
   \put(4,40.9){$0$}   \put(4,78){$\frac{1}{5}$}
   \put(15.3,1){$-\frac{1}{2}$}   \put(49.75,2){$0$}   \put(75,1){$\frac{1}{2}$}
   \put(77,87){$\pi_4$}
 \end{picture}
\]
\caption{Conics in the planes $\pi_i$}
 \label{F:conics}
\end{figure}

For each $i$, let $C_i$ be the conic $\pi_i\cap Q_i$, shown in the plane $\pi_i$
in Figure~\ref{F:conics}. 
Here, the horizontal coordinate is $t$, the parameter of the line $\ell_i$, while the
vertical coordinate is $y{-}1$ for $\pi_1$,  $y{+}1$ for $\pi_2$, $z{-}\frac{1}{4}$ for
$\pi_3$, and $z$ for $\pi_4$.

For each $i=1,\dotsc,4$, rotate line $\ell_i$ in plane $\pi_i$ very slightly
about a point that is far from the conic $C_i$, obtaining a new line $k_i$ in
$\pi_i$ which also meets $C_i$ in two points. 
Consider now the transversals to $\ell_i\cup k_i$, for $i=1,\ldots,4$.
Because $k_i$ is near to $\ell_i$ and there were two transversals to 
$\ell_1,\ell_2,\ell_3,\ell_4$, there will be 2 transversals to each of the 16 quadruples
of lines obtained by choosing one of $\ell_i$ or $k_i$ for $i=1,\ldots,4$.
By our choice of the point of rotation, all of these will meet $\ell_i$ and
$k_i$ in one of the two thin wedges they form.
In this wedge, form a triangle by adding a third side so that the edges on 
$\ell_i$ and $k_i$ contain all the points where the transversals meet the lines.
The resulting triangles will then have at least 32 common tangents.
We claim that by carefully choosing the third side (and tuning the rotations) we are
able to get 30 additional tangents. 

\begin{figure}[htb]
\[
  \setlength{\unitlength}{1pt}
  \begin{picture}(160,150)(-10,0)
    \put(0,10){\includegraphics[width=5cm]{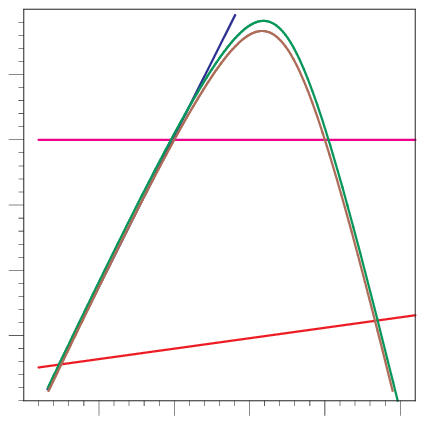}}
    \put(27,0){$-1$} \put(52,0){$-\frac{1}{2}$}
    \put(86,0){$0$}  \put(110,0){$\frac{1}{2}$}
    \put(135,0){$1$} 
    \put(-12,35){$-\frac{3}{10}$} \put(-12,56){$-\frac{2}{10}$}
    \put(-12,77){$-\frac{1}{10}$}  \put(-2,98){$0$}
    \put(-4,119.5){$\frac{1}{10}$}  
    \put(86,43){$k_4$} \put(128,106){$\ell_4$}
    \put(45,135){$m_4$}  \put(63,137){\vector(1,0){15}}
   \end{picture}
  \qquad
  \begin{picture}(160,150)(-10,0)
    \put(0,10){\includegraphics[width=5cm]{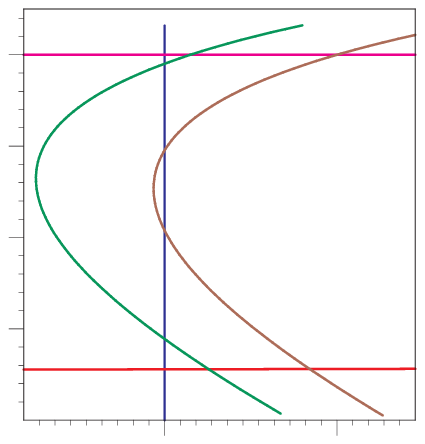}}  
    \put(40,-1){$-\frac{511}{1000}$}  \put(99,-1){$-\frac{1}{2}$} 
    \put(-12,42){$-\frac{3}{10}$} \put(-12,70){$-\frac{2}{10}$}
    \put(-12,98.4){$-\frac{1}{10}$}  \put(-2,126){$0$}
    \put(123,39){$k_4$} \put(123,120){$\ell_4$}
    \put(75,85){$m_4$}  \put(73,88){\vector(-1,0){13}}
    \put(95,67){$C_4$}  \put(93,66){\vector(-1,-1){13}}
 \end{picture}
\]
\caption{Configuration in plane $\pi_4$}
 \label{F:Pi_4}
\end{figure}

To begin, look at Figure~\ref{F:Pi_4} which displays the configuration in $\pi_4$ given
by the four triangles from Table~\ref{TA:triangles}.
Since the lines $\ell_i$ and $k_i$ for $i=1,2$ are extremely close, 
the four conics given by transversals to them and to $\ell_3$ cannot be resolved in
these pictures.
The same is true for the four conics given by $k_3$, so that each of the
apparent 2 conics are clusters of four nearby conics. 
The picture on the left is a view of this configuration in the coordinates for
$\pi_4$ of Figure~\ref{F:conics}.
It includes a secant line $m_4$ to the conics.
We choose coordinates on the right so that $m_4$ is vertical, but do not change the
coordinates on $\ell_4$.
The horizontal scale has been accentuated to separate the two clusters of conics.
The three lines, $\ell_4$, $k_4$, and $m_4$ form the triangle $t_4$.
Let its respective edges be $e_4$, $f_4$, and $g_4$.
Each edge meets each of the 8 conics in two points and these 48 points of
intersection give 48 lines tangent to the four triangles.

This last assertion that the 16 lines transversal to $m_4$ and to 
$\ell_i\cup k_i$ for $i=1,2,3$ meet the edges of the triangles $t_1$, $t_2$, and $t_3$
needs justification.
Consider for example the transversals to $\ell_1$, $\ell_2$, and $\ell_3$.
These form a ruling of the doubly-ruled quadric $Q_4$ and are parameterized by their
point of intersection with $\ell_1$.
The intersection of $Q_4$ with $\pi_4$ is the conic $C_4$.
Since the intersections of the conic $C_4$ with the segment $g_4$ supported on
$m_4$ lie between its intersections with $\ell_4$ and $k_4$, the corresponding
transversals to  $\ell_1$, $\ell_2$, $\ell_3$, and $g_4$ meet $\ell_1$ between
points of $\ell_1$ met by common transversals to $\ell_4\cup k_4$ and $\ell_1$,
$\ell_2$, and $\ell_3$. 
The same argument for the other lines and for all 8 conics justifies the assertion.

Na{\"\i}vely, we would expect that this same construction (the third side cutting all 8
conics in $\pi_i$) could work to select each of the remaining sides of the triangles
$g_3$, $g_2$, and $g_1$, and that this would give four triangles having 
\[
   32+16+16+16+16\ =\ 96
\]
common tangents. 
Unfortunately this is not the case.
In the earlier conference version of this paper, we gave a construction that we
claimed would yield 88 common tangents.
Attempting that construction using Maple revealed a flaw in the
argument and the current construction of four triangles with 62 common tangents
is the best we can accomplish.

In $\pi_4$, the conics come in two clusters, depending upon whether or not
they correspond to $\ell_3$ or to $k_3$.
In order for the edge $g_4$ to cut all conics, the angle between $\ell_4$ and $k_4$ has
to be large, in fact significantly larger than the angle between $\ell_3$ and $k_3$.
Thus in $\pi_3$, the conics corresponding to $\ell_4$ are quite far from the conics
corresponding to $k_4$, and the side $g_3$ can only be drawn to cut four of the conics,
giving 8 additional common tangents.
Similarly, $g_2$ can only cut two conics, and $g_1$ only 1.
In this way, we arrive at four triangles having
\[
  32+16+8+4+2\ =\ 62
\]
common tangents, which we can verify by computer.
\QED

\section{Upper bound for triangles in $T$}\label{sec:upper}
Four triangles in $T$ have at most 162 common tangents.
If the triangles are disjoint, we slightly improve this upper bound to 
156.
Our method will be to show that not all $81=3^4$ quadruples of edges 
can give rise to a common tangent.
Our proof follows that for the upper bound on the number of
tangents to four polytopes~\cite{cccg02}, limiting the number of
configurations for disjoint triangles in $\mathbb{R}^3$. 
We divide the proof into two lemmas, which do not assume that the triangles lie in $T$. 

In order for a tangent to meet an edge $e$, the plane it spans with $e$ must meet one
edge from each of the other triangles.
A triple of edges, one from each of the other triangles, is \Blue{\emph{contributing}}
if there is a plane containing $e$ which meets the three edges.
We say that an edge $e$ \Blue{\emph{stabs}} a triangle $t$ if its supporting
line meets the interior of $t$.

\begin{lemma}\label{L:26}
 Let $e$ be an edge of some triangle.
 If $e$ stabs exactly one of the other triangles, then there are at most $26$
 contributing triples of edges.
 If $e$ stabs no other triangle, then there are at most $25$
 contributing triples.
\end{lemma}

It is not hard to see that if $e$ stabs at least two of the other triangles, then each
of the $27=3^3$ triples of edges can be contributing.\medskip

\begin{figure}[htb]
 \[
   \begin{picture}(150,56)(0,-20)
    \put(0,0){\includegraphics[width=5cm]{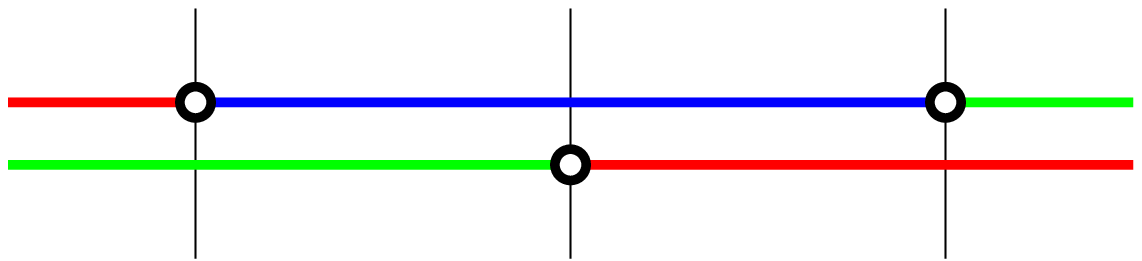}}
    \put(8,27){\Red{3}} \put(45,27){\Blue{1}} 
       \put(90,27){\Blue{1}} \put(125,27){\Green{2}} 
    \put(8,-2){\Green{2}} \put(45,-2){\Green{2}} 
       \put(90,-2){\Red{3}} \put(125,-2){\Red{3}}
       \put(50,-20){$e$ stabs $t$}
   \end{picture}
   \qquad\qquad
   \begin{picture}(150,56)(0,-22)
    \put(0,0){\includegraphics[width=5cm]{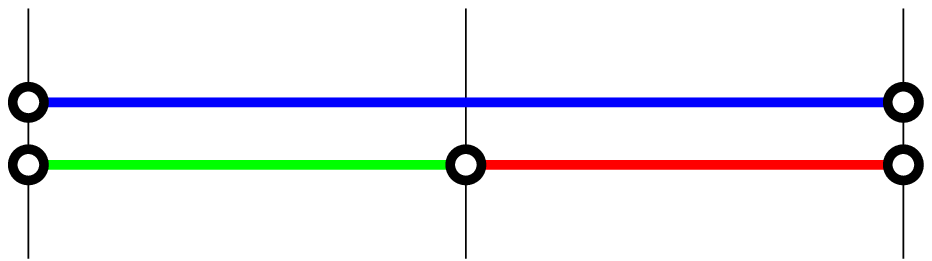}}
    \put(40,27){\Blue{1}}  \put(95,27){\Blue{1}}
    \put(40,-2){\Green{2}}  \put(95,-2){\Red{3}}
       \put(27,-20){$e$ does not stab $t$}
   \end{picture}
 \]
 \caption{Stabbing and non-stabbing configurations}
 \label{F:stab}
\end{figure}

\noindent{\it Proof.}
 Suppose that $e$ is an edge of some triangle.
 Let $\pi(\alpha)$ be the pencil of planes containing $e$.
 (This is parametrized by the angle $\alpha$.) 
 For each edge $f$ of another triangle $t$, there is an interval of angles $\alpha$ for
 which $\pi(\alpha)$ meets $f$.
 Figure~\ref{F:stab} illustrates the two possible configurations for these intervals,
 which depend upon whether or not $e$ stabs the triangle $t$.
 The intervals are labeled 1, 2, and 3 for the three edges of $t$.
 When $e$ stabs $t$, these intervals cover the entire range of $\alpha$ and the picture
 is actually wrapped.
 Call this a \Blue{{\it stabbing diagram}}.
 When the supporting line of $e$ does not meet $t$, these intervals do not cover the
 entire range of $\alpha$, and there are two endpoints and one \Blue{{\it interior
 vertex}} of the diagram.
 If the supporting line of $e$ meets an edge of $t$, then the two endpoints of the
 non-stabbing diagram wrap around and coincide.
 Call either of these last two configurations a non-stabbing diagram.

\begin{figure}[htb]
\[
  \begin{picture}(150,130)(0,-18)
   \put(0,10){\includegraphics[width=5cm]{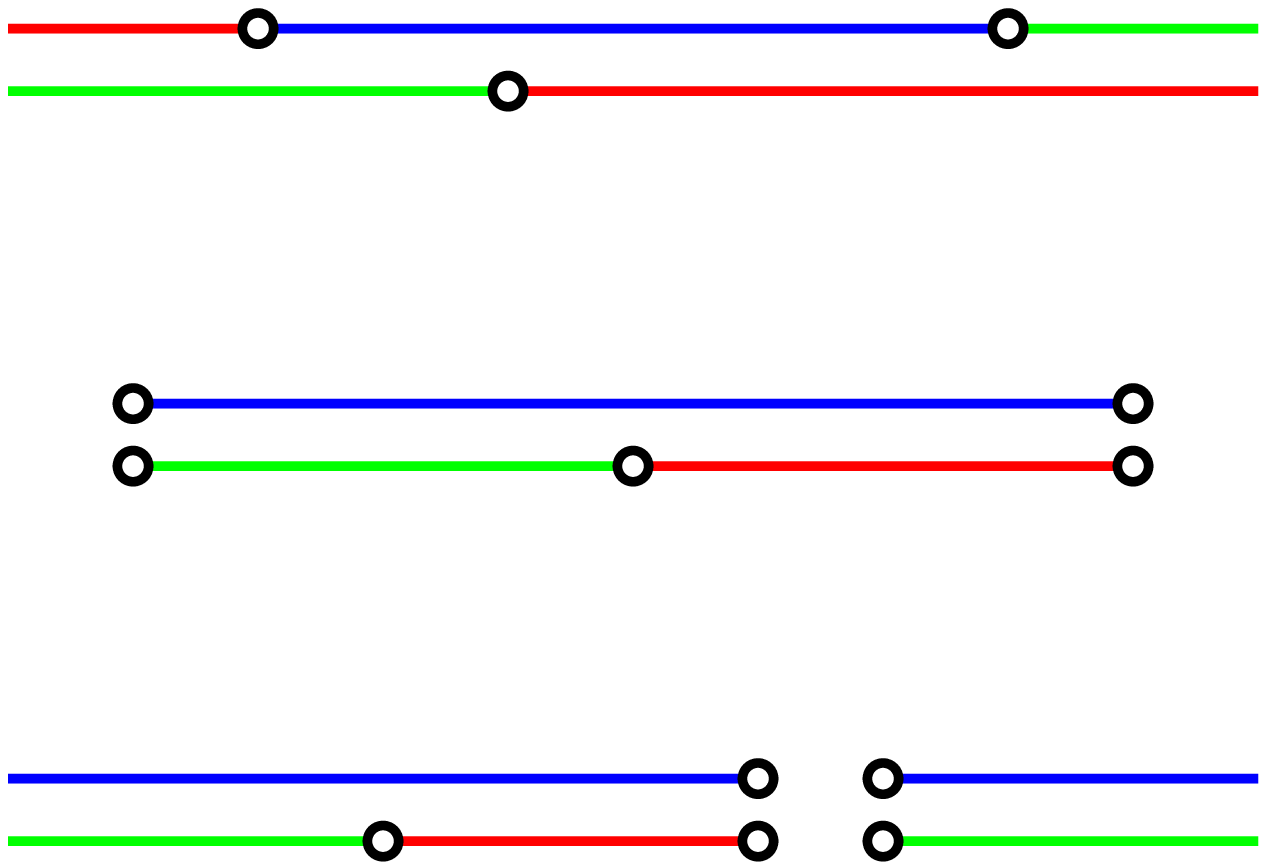}}
    \put(10,109){\Red{3}}  \put(79,109){\Blue{1}}  \put(125,109){\Green{2}} 
   \put(25,84){\Green{2}} \put(100,84){\Red{3}}
          \put(68,67){\Blue{1}}
   \put(40,42){\Green{2}} \put(95,42){\Red{3}}
   \put(45,25){\Blue{1}}  \put(120,25){\Blue{1}}
   \put(20,0){\Green{2}} \put(60,0){\Red{3}} \put(120,0){\Green{2}} 
   \put(30,-18){$e$ stabs one triangle}
   \end{picture}
    \qquad\quad
  \begin{picture}(150,130)(0,-18)
   \put(0,10){\includegraphics[width=5cm]{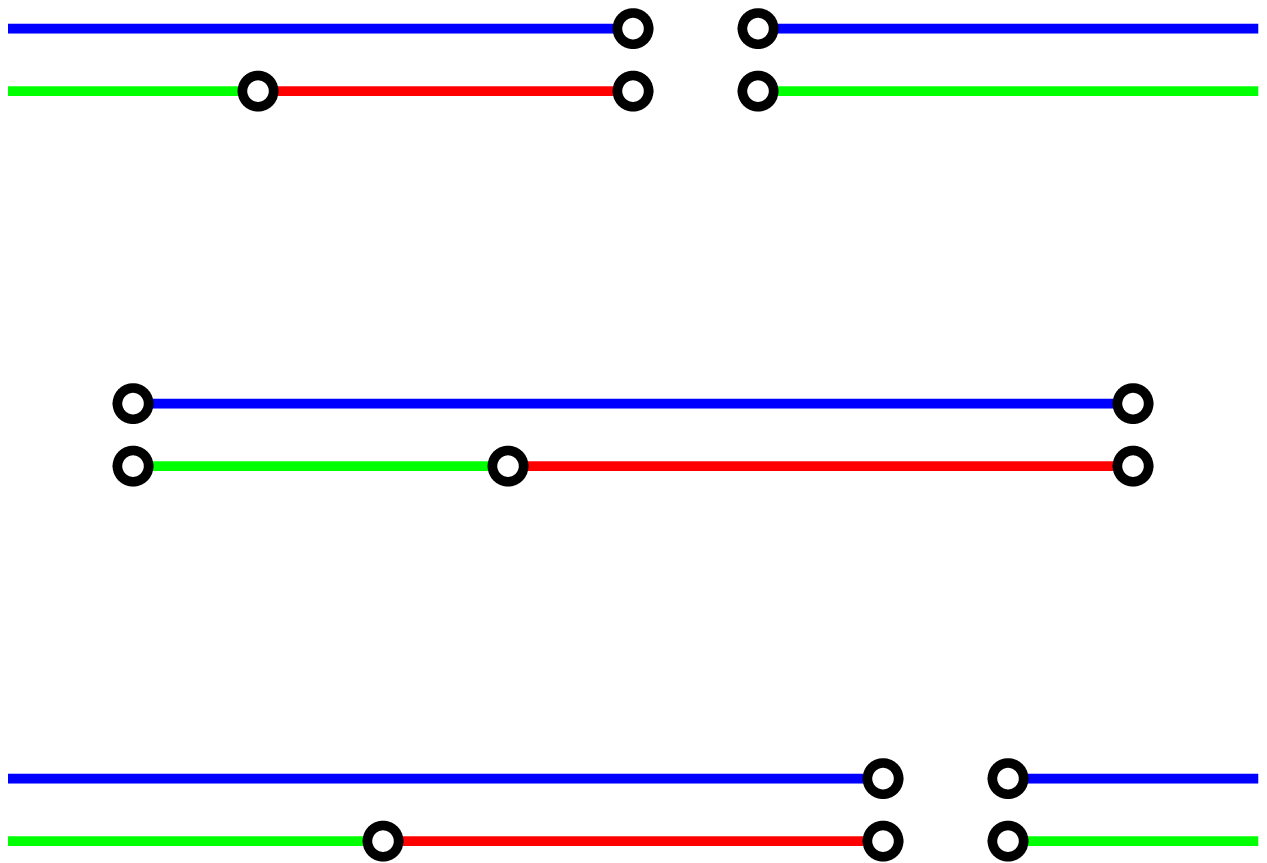}}
   \put(35,109){\Blue{1}}  \put(110,109){\Blue{1}}
   \put(10,84){\Green{2}} \put(47,84){\Red{3}} \put(110,84){\Green{2}} 
          \put(70,67){\Blue{1}}
   \put(35,42){\Green{2}} \put(90,42){\Red{3}}
   \put(60,25){\Blue{1}}  \put(125,25){\Blue{1}}
   \put(20,0){\Green{2}} \put(70,0){\Red{3}} \put(125,0){\Green{2}} 
   \put(35,-18){$e$ stabs no triangle}
   \end{picture}
\]
 \caption{Configurations with 26 and 25 contributing triples}
 \label{F:26}
\end{figure}

 To count contributing triples, we line up (overlay) diagrams from each of the three
 triangles not containing $e$  and count how many of the 27 triples
 $\{1,2,3\}^3$, one from each triangle, occur at some value of $\alpha$.
 For example,  Figure~\ref{F:26} displays a configuration with 26 contributing triples
 (where $e$ stabs a single triangle) and a configuration with 25 contributing triples
 ($e$ stabs no other triangles).
 The configuration on the left is missing the triple $(2,3,3)$, while the configuration on
 the right is missing the triples $(2,2,3)$ and $(3,3,2)$.

 These configurations are the best possible.
 Indeed, begin with two non-stabbing diagrams in which all 
 9 pairs of edges occur.
 (If only 8 pairs occurred, there would be at most 24 contributing triples.)
 The unique way to do this up to relabeling the edges is given by the lower
 two diagrams in either picture in Figure~\ref{F:26}. 
 These two diagrams divide the domain of $\alpha$ into 6 intervals (the two at the ends
 are wrapped). 
 The five pairs involving 1 occur in two intervals, but four exceptional pairs
 $\{(2,2),(2,3),(3,2),(3,3)\}$ occur uniquely in different intervals.

 Consider now a third diagram.
 An exceptional pair extends to three contributing triples
 only if all three sides in the third diagram meet the interval corresponding to that
 pair. 
 If the third diagram is stabbing, then one of its three vertices lies in
 that interval---thus there is at least one triple which does not contribute.
 If the third diagram is non-stabbing, then either the middle vertex or else both
 endpoints must lie in that interval---thus there are at least two triples which
 do not contribute.
\QED

\begin{lemma}
 At most $78$ quadruples of edges of four disjoint triangles can lead to a common
 tangent. 
\end{lemma}

\begin{figure}[htb]
\[
  \begin{picture}(150,110)
   \put(0,0){\includegraphics[width=5cm]{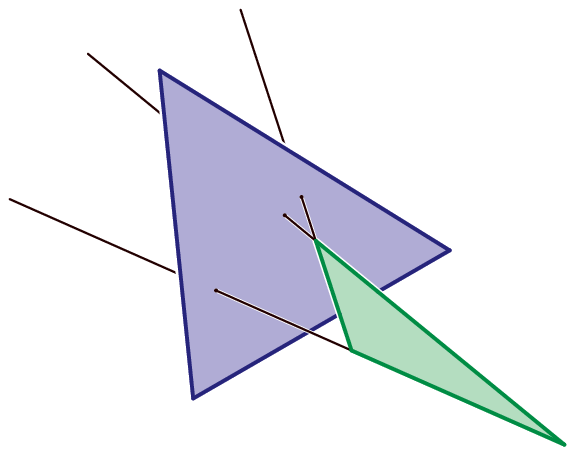}}
   \put(20,60){$t'$}   \put(130,20){$t$}
    \put(65,0){(a)}
  \end{picture}
   \qquad
  \begin{picture}(150,110)
   \put(0,-0.7){\includegraphics[width=5.5cm]{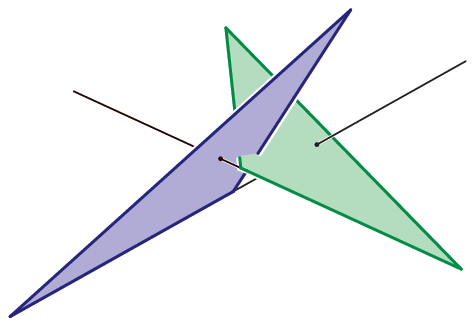}}
    \put(65,0){(b)}
  \end{picture}
\]
 \caption{(a) Two disjoint triangles can have at most 3 stabbing lines.\newline
          \mbox{\ }\hspace{2.05cm}(b) Two intersecting triangles may have up to four.}
 \label{F:stabbing}
\end{figure}

\noindent{\it Proof.}
 First consider the maximum number of stabbing edges between two triangles.
 If the triangles are disjoint, then there are at most three stabbing edges; 
 one triangle could have three edges stabbing the other.
 Indeed, if at least two supporting lines of a triangle $t$ meet another triangle
 $t'$ which is disjoint from $t$, then $t$ lies entirely on one side of the plane
 supporting $t'$, and thus no supporting lines of $t'$ can meet $t$.
 Figure~\ref{F:stabbing}(a) shows a configuration in which all three supporting lines
 of $t$ stab $t'$. 

 Consider now the bipartite graph between 12 nodes representing the edges of the four
 triangles and 4 nodes representing the triangles.
 This graph has an arc between an edge $e$ and a triangle $t$ if the line supporting $e$
 stabs $t$. 
 (We assume that $e$ is not an edge of $t$.)
 We just showed that the edges of one triangle $t$ can have at most 3 arcs incident on
 another triangle $t'$, and so this graph has at most 18 edges.
 
 Let the weight of a triangle be the number of arcs emanating from its edges in this
 graph.
 As the graph has at most 18 arcs, at least one triangle has weight less than 5.
 We argue that there is a triangle of weight at most 3.
 This is immediate if the graph has 15 or fewer edges.
 On the other hand, this graph has more structure.
 If it has 18 edges, then all pairs of triangles are in the configuration of
 Figure~\ref{F:stabbing}(a), and so every triangle has weight a multiple of 3, which
 implies that some triangle has weight at most 3.
 If the graph has 17 edges, then there is exactly one pair of triangles with only two
 stabbing edges, and so the possible weights less than 5 are 0, 2, and 3.
 If the graph has 16 edges, then there is one pair with only one edge stabbing,
 or two pairs with 2 edges stabbing.
 There can be at most 2 triangles of weight 4, and again we conclude that there is
 triangle with weight at most 3.

 If a triangle has weight at most three, either all three edges stab a unique triangle,
 or else one edge stabs no triangles and another edge stabs at most one other triangle.
 We sum the number of contributing triples over the edges of this triangle.
 By Lemma~\ref{L:26}, this sum will be at most 26+26+26=78 if all three edges stab a
 unique triangle and at most 27+26+25=78 if not.
 This proves the lemma.
\QED

\begin{remark}
 There exist four disjoint triangles whose bipartite graph has exactly 18 edges.
 Thus the previous argument cannot be improved without additional ideas. 
 It is conceivable that further restrictions the bipartite graph may exist, 
 leading to a smaller upper bound.
\end{remark}

\begin{remark}
 This proof does not enable us to improve the bound when the triangles are not disjoint.
 Two intersecting triangles can induce up to
 four arcs (see Figure~\ref{F:stabbing}(b)) and thus the total number of arcs is bounded
 above by 24.
 The minimal weight of a triangle is then 6, and the edges of such a triangle could all
 have degree 2, which leads to no restrictions.
\end{remark}

\section{Upper bound on the number of components}\label{S:upper-comp}

Let $\mathcal{F}$ and $\mathcal{I}$ be the sets of quadruples of edges, one from
each of four triangles, whose supporting lines have finitely and infinitely,
respectively, many common transversals.  Let $n_\mathcal{F}$ and $n_\mathcal{I}$
be the sum over all quadruples of edges in $\mathcal{F}$ and $\mathcal{I}$,
respectively, of the numbers of connected components of common transversals to
each quadruple of edges.  Note that the number of quadruples in $\mathcal{F}$ and
$\mathcal{I}$ is $|\mathcal{F}| + |\mathcal{I}|=81$.

Consider a connected component $c$ of common transversals to a quadruple of
edges $q\in\mathcal{I}$. The arguments of~\cite{BELSW} show that $c$
contains a line that meets a vertex of one of the four edges. That line is thus
transversal to another quadruple $q'$ of edges.  
Thus, the connected component $c$ of common transversals to $q$ is connected with
a connected component $c'$ of common transversals to $q'$. 
If $q'\in\mathcal{F}$ we charge the component $c\cup c'$ to $c'$.  Otherwise $q$
and $q'$ are both in $\mathcal{I}$ and the component $c\cup c'$ is counted twice.
The number of connected components of tangents to  four triangles is
thus at most $n_\mathcal{F} + {n_\mathcal{I}}/{2}$.

Since any four lines admit at most 2 or infinitely many
transversals, $n_\mathcal{F}\leq 2|\mathcal{F}|$. Also,   any four segments
admit at most 4 connected components of common transversals~\cite{BELSW}, thus 
$n_\mathcal{I}\leq 4|\mathcal{I}|$. 
Hence, the number of connected components of tangents to four triangles is 
at most $2|\mathcal{F}| + 2|\mathcal{I}| = 162$.

This still may overcount the number of connected components of
tangents, but further analysis is very delicate.
Such complicated arguments are not warranted as we have already obtained the
upper bound of 162 common tangents to four triangles in $T$.
As in Section~\ref{sec:upper}, if the triangles are disjoint, then not all
quadruples of edges can contribute, which lowers this bound to 156.

\section{Random triangles}\label{S:search}
 We proved Theorem~\ref{T:lower} by exhibiting four triangles having 62 common
 tangents.
 We do not know if that is the best possible.
 Since the geometric problem of determining the tangents to four triangles
 is computationally feasible---it is the disjunction of 81 problems with
 algebraic degree 2 and simple inequalities on the solutions---we investigated
 it experimentally.

 For this, we generated 5\,000\,000 quadruples of triangles whose
 vertices were points with integral coordinates chosen uniformly at random from the cube
 $[-1000,1000]^3$.
 For each, we computed the number of tangents.
 The resulting frequencies are  recorded in Table~\ref{TA:random}.
\begin{table}[htb] 

\[
    \begin{tabular}{|l||c|c|c|c|c|c|c|c|}\hline
    Number    &  0    &  2   &  4   &  6   &   8  & 10   &  12     & 14\rule{0pt}{13pt}\\\hline
    Frequency &1\,515\,706&331\,443&646\,150&403\,679&637\,202&327\,159&358\,312&238\,913\rule{0pt}{13pt}\\\hline
   \end{tabular}\vspace{.25cm}  
\]
\[
   \begin{tabular}{|c|c|c|c|c|c|c|c|c|c|c|c|c|}\hline
     16   &  18  & 20  & 22  & 24  & 26  & 28 & 30 &32 &34 &36&38&40\rule{0pt}{13pt}\\\hline
    253\,396&114\,046&80\,199&44\,870&27\,726&12\,426&5\,796&2\,016&813&111&30&3&4\rule{0pt}{13pt}\\\hline
   \end{tabular}\vspace{.25cm}  
\]
 \caption{Number of triangles with a given number of tangents, out of 5\,000\,000 randomly constructed triangles\label{TA:random}}

\end{table}
 This search consumed 17 million seconds of CPU time on 1.2GHz processors at the
 MSRI and a DEC Alpha machine at the University of Massachusetts in 2004. 
 It is archived on the web 
 page\footnote{{\tt www.math.tamu.edu\~{}sottile/stories/4triangles/index.html}}
 accompanying this article.

\begin{figure}[htb] 
\[
   \includegraphics[width=5.5cm]{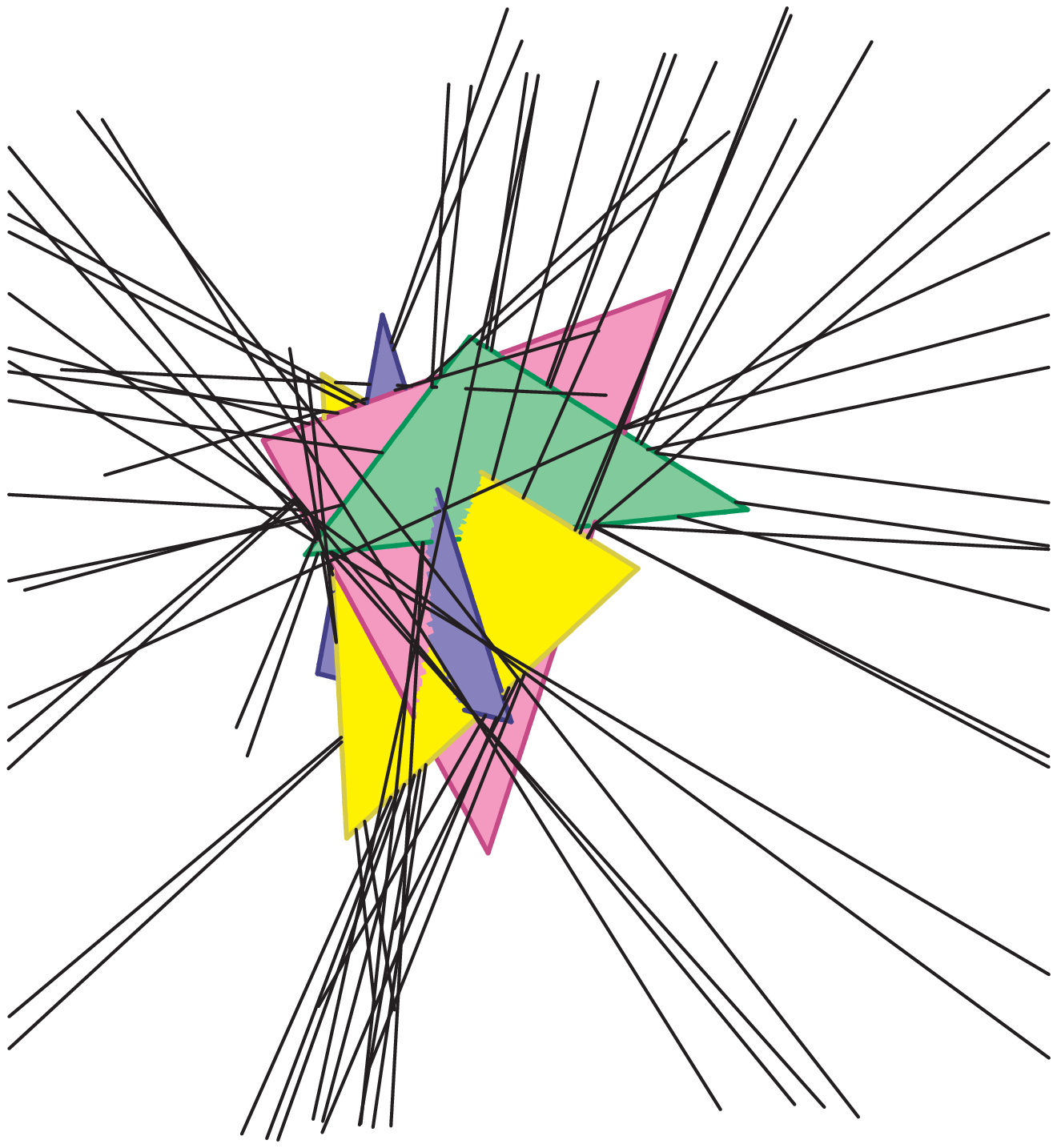}\qquad
   \begin{picture}(160,180)
    \put(0,0){\includegraphics[width=5.5cm]{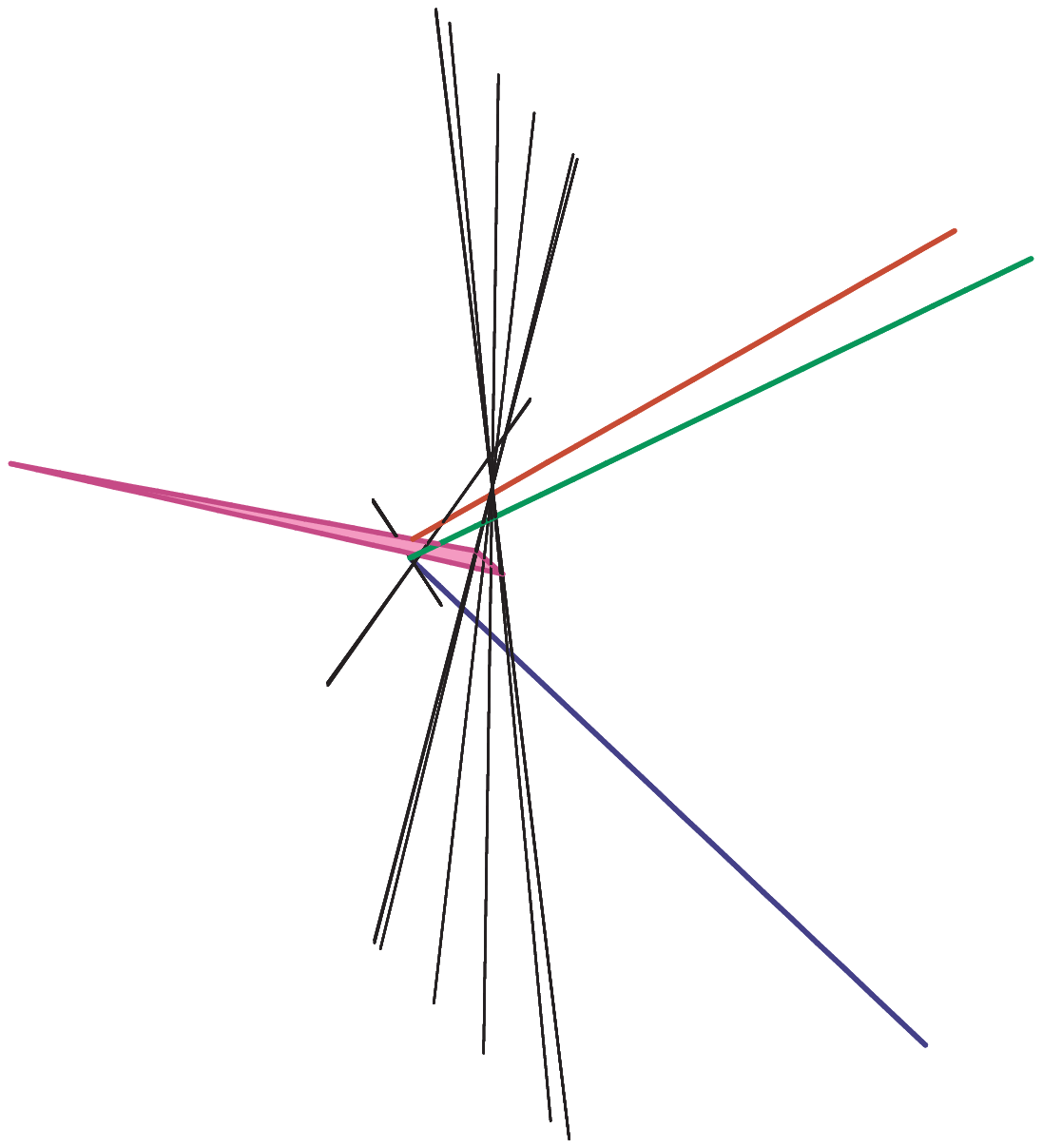}}
    \put(130, 30){$t_1$}
    \put(122,138){$t_2$}
    \put(140,115){$t_3$}
    \put(  5,108){$t_4$}
   \end{picture}
\]
\caption{Triangles with many common tangents\label{F:40}}
\end{figure}

 In this search, we found four different quadruples of triangles with 40 common
 tangents, and none with more.
 The vertices of one are given in Table~\ref{TA:40}.
\begin{table}[htb] 
 \begin{center}
    \begin{tabular}{|c||c|c|c|}\hline
   Triangle&\multicolumn{3}{c|}{Vertices}\\\hline
   $t_1$&$(-4, -731, -336) $&$(297, -507, 978)$&$( 824, -62, -359)$\\\hline
   $t_2$&$(531, -631, -820)$&$(-24, -716, 713)$&$( 807, 377,  177)$\\\hline
   $t_3$&$(586, -205,  952)$&$(861, -774, 235)$&$(-450, 758,  161)$\\\hline
   $t_4$&$(330, -141, -908)$&$(942, -920, 651)$&$(-226, 489,  968)$\\\hline
  \end{tabular}
 \end{center}
 \caption{Four triangles with 40 common tangents\label{TA:40}}

\end{table}
These triangles are rather `fat', in that none have very small angles.
Contrast that to the triangles of our construction in Section~\ref{sec:upper}.
In Figure~\ref{F:40} we compare these two configurations of triangles.
On the left is the configuration of triangles from Table~\ref{TA:40}, 
together with their 40 common tangents, while on the right is the configuration of
triangles having 62 common tangents.
The triangles are labeled in the second diagram, as they are hard to distinguish
from the lines.
As we remarked in Section~\ref{sec:upper}, many of the lines are extremely close
and cannot be easily distinguished; that is why one can only count 8 lines in
this picture.

\section*{Acknowledgments}
  This research was initiated at the Second McGill-INRIA Workshop on Computational
  Geometry in Computer Graphics, February 7--14, 2003,  co-organized by
  H. Everett, S. Lazard, and S. Whitesides, and held at the Bellairs
  Research Institute of McGill University.
We would like to thank the other participants of the workshop for useful discussions.


\bibliographystyle{abbrv}

\def\baselinestretch{1.0}
\normalsize

\end{document}